\documentclass{ifacconf}

\usepackage{graphicx}      
\usepackage{natbib}        
\usepackage{amsmath,amssymb,amsfonts}
\usepackage{enumitem}
\usepackage{xcolor}

\newtheorem{theorem}{Theorem}

\newtheorem{definition}{Definition}
\newtheorem{assumption}{Assumption}

\newtheorem{remark}{Remark}

\newcommand{\R}{\mathbb{R}}

\newcommand{\bref}[1]{(\ref{#1})}
\newcommand{\tref}[1]{Theorem~\ref{#1}}
\newcommand{\lref}[1]{Lemma~\ref{#1}}
\newcommand{\aref}[1]{Assumption~\ref{#1}}
\newcommand{\dref}[1]{Definition~\ref{#1}}
\newcommand{\sref}[1]{Section~\ref{#1}}
\newcommand{\rref}[1]{Remark~\ref{#1}}
\newcommand{\fref}[1]{Figure~\ref{#1}}

\begin{document}
\begin{frontmatter}

\title{SOS Construction of Compatible Control Lyapunov and Barrier Functions}

\author[First]{Michael Schneeberger} 
\author[Second]{Florian D{\"o}rfler}
\author[First]{Silvia Mastellone} 

\address[First]{Institute of Electrical Engineering, FHNW, Windisch, Switzerland (e-mails: michael.schneeberger@fhnw.ch, silvia.mastellone@fhnw.ch)}
\address[Second]{Department of Information Technology and Electrical Engineering, ETH Z{\"u}rich, Z{\"u}rich, Switzerland, (e-mail: dorfler@control.ee.ethz.ch)}

\begin{abstract}                

   We propose a novel approach to certify closed-loop stability and safety of a constrained polynomial system based on the combination of Control Lyapunov Functions (CLFs) and Control Barrier Functions (CBFs).
   For polynomial systems that are affine in the control input, both classes of functions can be constructed via Sum Of Squares (SOS) programming.
   Using two versions of the Positivstellensatz we derive an SOS formulation seeking a rational controller that --- if feasible --- results in compatible CLF and multiple CBFs.
\end{abstract}

\end{frontmatter}

\section{Introduction}

When dealing with systems that have state constraints, it is crucial to have a controller that ensures both stability and compliance with the constraints.
In most cases, feedback control design focuses mainly on achieving stability, while protection functions only engage when constraints are violated.
Unfortunately, this approach results in downtime and the need to investigate faults.
By contrast, a controller that is both stable and safe can trade off control performance for the ability to prevent unsafe states.

The focus of this paper is on constructing compatible CLFs and CBFs that can certify both stability and safety in control systems.
A CLF establishes conditions for the existence of a stabilizing controller for a given control system.
Similarly, a CBF guarantees the existence of a controller that can render the control system safe.
According to \cite{WielandAllgower:07}, a system is considered safe if any state trajectory starting from a safe set of states remains within an allowable region defined by the state constraints.
For systems that are affine in the input, a controller that meets the CLF and CBF conditions can be implemented by solving an online Quadratic Programming (QP), see \cite{AmesCoogan:19}.
However, for some states, these conditions may conflict with each other.
For both to hold jointly for all states, the control-sharing property is additionally required, see \cite{GramBlanch:13,Xu:16}.
Finding such compatible CLF and CBF is generally difficult.
For polynomial systems, however, this can be achieved by formulating SOS constraints and solving them via Semidefinite Programming (SDP).

The contributions of this paper is twofold:
first, we derive SOS constraints on compatible CLF and multiple CBFs using two versions of the Positivstellensatz, and second, an algorithm is developed that efficiently finds solutions to these SOS constraints by maximizing a surrogate of the volume of the safe set.
The conditions on compatible CBFs are formulated in \cite{IsalyGhanbarpour:22, TanDimarogonas:22} but without giving a method to construct them.
A constructive approach described by \cite{Clark:21} is based on the introduction of additional SOS constraints to enforce compatibility.
In this paper, we reveal a correspondence between the SOS constraint derived from the CLF, resp. CBF, condition, and the existence of a rational controller that renders the closed-loop system stable, resp. safe.
By restricting such controllers to be identical, we derive a new set of SOS constraints that guarantee compatibility between a CLF and multiple CBFs without the introduction of additional SOS constraints.

These SOS constraints contain bilinear terms, and hence, cannot be directly converted to an SDP.
We therefore present an alternating algorithm that searches simultaneously for a CLF and multiple CBFs by repeatedly solving two SDPs.
In particular, the algorithm seeks to maximize the volume of the safe set, which is given by the intersection of the invariant sets defined by each CBF.
Multiple CBFs offer additional flexibility to increase the volume of the safe set when a single CBF does not suffice.
Similar approaches were presented in \cite{AnghelMilano:13,KunduGeng:19,WangMargellos:22} but they either only searched for a single CLF or a single CBF.
\cite{KordaHenrion:14} proposed another intriguing approach for identifying a safe set using an infinite-dimensional linear programming problem.
We demonstrate the utility of our approach with a power converter control example for which safety is of paramount importance.

The paper is structured as follows: 
In \sref{sec:preliminaries}, we introduce the notation adopted in the paper and review some preliminaries.
Then, we recall the definition of a CLF and CBF and derive a rational controller in \sref{sec:stability_safety}.
In \sref{sec:unification}, we combine the CLF and CBF by introducing the control-sharing property.
\sref{sec:sos_program} defines the SOS program encoding the CLF and CBF conditions, and \sref{sec:algorithm} presents the algorithm that solves the SOS program.
Numerical simulation are given in \sref{sec:simulation}.
Finally, \sref{sec:conclusions} is dedicated to concluding remarks and future work.

\section{Preliminaries \& Notation} \label{sec:preliminaries}

\subsection{Notation}

Across this paper, we adopt the following notation.
The shorthand $[t] := \lbrace1, 2, ..., t\rbrace$ is used to denote a range of numbers.
A scalar function $V: \R^n \to \R$ is positive definite w.r.t. $x^*$, if $V(x^*) = 0$ and $V(x) > 0$ for $x \neq x^*$.
$R[x]$ refers to the set of scalar polynomials in variables $x \in \R^n$, and $\Sigma[x]$ refers to the set of scalar SOS polynomials in $x$. 
A polynomial $p(x)$ is an SOS polynomial if it can be written as $p(x) = \sum_{i=1}^k g_i(x)^2$ for $g_i(x) \in R[x]$ and $i \in [k]$.
If $p(x) \in \Sigma[x]$, it can be expanded to
\begin{align*}
   p(x) = Z(x)^T Q Z(x),
\end{align*}
where $Z(x)$ is a vector of monomials in $x$, and $Q$ is a square positive semidefinite matrix.
In the following, we consider a polynomial control system described as
\begin{align} \label{eq:poly_system}
   \dot x = \mathcal F_c(x, u) = f(x) + G(x) u_c,
\end{align}
where $f(x) \in \left( R[x] \right)^n$ and $G(x) \in \left( R[x] \right)^{n \times m}$ are polynomial matrices, and $u_c \in \R^m$ is the control input vector.
System~\bref{eq:poly_system} with polynomial state feedback control policy $u(x) \in \left( R[x] \right)^m$ results in a closed-loop polynomial system of the form
\begin{align} \label{eq:autonomous_poly_system}
   \dot x = \mathcal F_a(x) = f(x) + G(x) u(x)
\end{align}
The state $x^*$ is called an equilibrium of~\bref{eq:autonomous_poly_system}, if $\mathcal F_a(x^*) = 0$.

A set $\mathcal X$ is called forward invariant \cite[]{Khalil:02} with respect to \bref{eq:autonomous_poly_system} if for every $x_0 \in \mathcal X$, $x(t) \in \mathcal X$ for all $t \in \R_+$.
We say that a system \bref{eq:autonomous_poly_system} is safe \cite[]{WielandAllgower:07} w.r.t. an allowable set of states $\mathcal X_a \subseteq \R^n$ and the safe states $\mathcal X_s \subseteq \R^n$, if $\mathcal X_s$ is forward invariant and $\mathcal X_s \subseteq \mathcal X_a$.

\subsection{SOS Programming}

An SOS program minimizing a quadratic cost function subject to SOS constraints is defined as follows:
\[
    \begin{array}{lll}
        \underset{v}{\mbox{minimize}} & v^T P v + c^T v \\
        \mbox{subject to} & p_i(x, v) \in \Sigma[x] & i=1,...,N.
    \end{array}
\]
Here $P \in \R^{l \times l}$ is a positive semi-definite matrix, $c \in \R^l$ is a vector, and $p_i(x, v) \in \Sigma[x]$ is an SOS polynomial in $x \in \R^n$ parameterized by $v \in \R^l$, i.e. it can be expanded as
\begin{align} \label{eq:param_poly}
   p_i(x, v) = Z_i(x)^T Q_i(v) Z_i(x),
\end{align}
where $Q_i(v) = Q_{0,i} + v_1 Q_{1,i} + ... + v_l Q_{1,l}$ is a square positive semidefinite matrix linearly parameterized by $v$ encoding the coefficients of the polynomial.

When deriving the SOS constraints in the following subsections, the parameter $v$ is omitted for simplicity, and the polynomial is simply written as $p_i(x)$ instead.

\subsection{Positivstellens{\"a}tze}

In the following, we present a version of the Positivstellensatz that can be regarded as a specialization of the weak Positivstellensatz. This version will be relevant for our analysis later on.
\begin{theorem} \label{thm:positivstellensatz}
   Given polynomials $f_1(x)$, $f_2(x)$, $g_1(x)$, $h_1(x)$, ..., $h_{m+1}(x) \in R[x]$ such that 
   \begin{align} \label{thm:positivstellensatz_emptyset}
      \begin{split}
          \Bigl\{ x \in \mathbb R^n \mid &f_1(x) \geq 0, f_2(x) \geq 0, g_1(x) \neq 0, \\
          &h_1(x) = 0, ..., h_{m+1}(x) = 0 \Bigr\} = \emptyset,
      \end{split}
   \end{align}
   then there exist polynomials \begin{itemize}[topsep=8pt,itemsep=5pt]
      \item $\begin{aligned}[t]
         f_\text{cone}(x) = s_1(x) f_1 (x) &+ s_2(x) f_2(x) \\
         &+ s_3(x) f_1(x) f_2(x)
      \end{aligned}$
      \item $h_\text{ideal}(x) = p_1(x) h_1(x) + ... + p_{m+1}(x) h_{m+1}(x)$
  \end{itemize} such that
  \begin{align} \label{thm:positivstellensatz_equality}
      -f_\text{cone}(x) - h_\text{ideal}(x) - g_1(x)^{2 k} \in \Sigma[x]
  \end{align}
  for all $x \in \mathbb R^n$, where $p_1(x), ..., p_{m+1}(x) \in R[x]$ are polynomials, $s_1(x), s_2(x), s_3(x) \in \Sigma[x]$ are SOS polynomials, and $k \in \mathbb N_+$.
\end{theorem}
\begin{pf}
   Using \cite[Theorem 4.4.2]{BochnakCoste:13}, we note that the cone $P$ generated by $f_1(x)$ and $f_2(x)$ is contained in $\Bigl\{ s_0(x) + s_1(x) f_1(x) + s_2(x) f_2(x) + s_3(x) f_1(x) f_2(x) \mid s_0, s_1, s_2, s_3 \in \Sigma[x] \Bigr\}$.
   \hfill$\blacksquare$
\end{pf}



If the set~\bref{thm:positivstellensatz_emptyset} is empty, the Positivstellensatz ensures that the polynomials $s_1(x), s_2(x), ..., p_{m+1}(x)$ exist without specifying their degree.
To find these polynomials computationally will require to iteratively increase the degree of the polynomials until a solution can be found for~\bref{thm:positivstellensatz_equality}. 
An increase in the degree of the polynomials will, however, deteriorate computation time.
For many practical examples of interest, however, low-degree polynomials suffice.

The following theorem states a version of the Positivstellensatz that can be seen as specializations of Putinar's Positivstellensatz.
\begin{theorem} \label{thm:putinar_positivstellensatz}
   Given polynomials $f_1(x)$, $f_2(x)$, $h_1(x)$, ..., $h_{m+1}(x) \in R[x]$ such that the set $\Bigl\{ x \in \R^n \mid f_2(x) \geq 0 \Bigr\}$ is compact and 
   \begin{align} \label{thm:putinar_positivstellensatz_emptyset}
      \begin{split}
          \Bigl\{ x \in \R^n \mid &f_1(x) \geq 0, f_2(x) \geq 0, \\
          &h_1(x) = 0, ..., h_{m+1}(x) = 0 \Bigr\} = \emptyset,
      \end{split}
   \end{align}
   then there exist polynomials \begin{itemize}[topsep=8pt,itemsep=5pt]
      \item $f_\text{cone}(x) = f_1 (x) + s_2(x) f_2(x)$
      \item $h_\text{ideal}(x) = p_1(x) h_1(x) + ... + p_{m+1}(x) h_{m+1}(x)$
   \end{itemize} such that
   \begin{align}
      \label{thm:putinar_positivstellensatz_equality}
      -f_\text{cone}(x) - h_\text{ideal}(x) \in \Sigma[x]
   \end{align}
   for all $x \in \R^n$, where $p_1(x), ..., p_{m+1}(x) \in R[x]$ are polynomials and $s_2(x) \in \Sigma[x]$ is an SOS polynomial.
\end{theorem}
\begin{pf}
   Let's consider the closed semialgebraic set 
   \begin{align*}
      \begin{split}
         K := \bigl\{ x \in \R^n \mid &f_2(x) \geq 0, \\
         &h_k(x) \geq 0 \quad k=1,...,m+1, \\
         -&h_k(x) \geq 0 \quad k=1,...,m+1 \bigr\},
      \end{split}
   \end{align*}
   and the quadratic module 
   \begin{align*}
      M(f_2(x), h_1(x), -h_1(x), ..., -h_{m+1}(x)),
   \end{align*}
   then $M$ is Archimedean according to \cite[Theorem 3.17]{Laurent:09} using the fact that $\bigl\{ x \in \R^n \mid f_2(x) \geq 0 \bigr\}$ is compact.
   The empty set condition \bref{thm:putinar_positivstellensatz_emptyset} is equivalent to the condition that $-f_1(x) > 0$ on K, and, therefore, $K$ is equivalent to
   \begin{align} \label{thm_proof:putinar}
      \begin{split}
         -f_1(x) = &s_0(x) + s_2(x) f_2(x) \\
         &+ \sum_{k=1}^{m+1} \left( s_{h,k}^p(x) - s_{h,k}^n(x) \right) h_k(x),
      \end{split}
   \end{align}
   where $s_0(x), s_2(x), s_{h,1}^p(x), ..., s_{h,m+1}^n(x) \in \Sigma[x]$ \cite[Theorem 3.20]{Laurent:09}.
   We note that $s_{h,k}^p(x) - s_{h,k}^n(x)$ can be replaced by a polynomial $p_k(x)$.

   \hfill$\blacksquare$
\end{pf}

\section{Closed-loop Stability and Safety} \label{sec:stability_safety}

In this section, we first review the concept of CLFs and CBFs.
Given a polynomial system~\bref{eq:poly_system}, we then propose an approach to construct a rational controller resulting from the SOS formulation of CLF and CBF conditions.
We prove global asymptotic stability and safety of the resulting closed-loop system~\bref{eq:autonomous_poly_system} given such a controller.

\subsection{Control Lyapunov Function}

For a control system~\bref{eq:poly_system} to be stabilizable around the equilibrium point $x^*$, we require the existence of a control input $u_c \in \R^m$ at every state $x \in \R^n$ that renders the sublevel sets of a scalar polynomial $V(x) \in R[x]$ forward invariant.
This motivates the following definition of a CLF:

\begin{definition} (\cite{Isidori:95}) \label{def:clf}
   Consider a differentiable function $V: \R^n \to \R$ such that $V(x^*) = 0$ and $V(x) > 0$ for all $x \neq x^*$.
   Such a scalar function is called a \underline{Control Lyapunov Function} (CLF) for the control system~\bref{eq:poly_system} if
   \begin{align} \label{eq:clf_condition}
      \nabla V(x)^T f(x) < 0
   \end{align}
   for all $x \in \Bigl\{ x \in \R^n \mid \nabla V(x)^T G(x) = 0, x \neq x^* \Bigr\}$.
\end{definition}
The inequality~\bref{eq:clf_condition} can also be formulated as the empty set condition
\begin{align} \label{eq:clf_emptyset}
   \begin{split}
       \Bigl\{ x \in \R^n \mid & \nabla V(x)^T f(x) \geq 0, \nabla V(x)^T G(x) = 0, \\
       &x \neq x^* \Bigr\} = \emptyset.
   \end{split}
\end{align}
Condition~\bref{eq:clf_emptyset} when restricted to \emph{polynomial} CLFs can be solved via SOS programming.
Hence, we restrict our attention to polynomial scalar functions $V \in R[x]$ for the rest of the paper.

Next, we replace the inequality constraints $x = (x_1, ..., x_n)\allowbreak \neq x^*$ in~\bref{eq:clf_emptyset} by a single inequality constraint $l(x) \neq 0$ (c.f. \cite{TanPackard:04}), where $l(x^*) = 0$ and $l(x) \neq 0$ elsewhere.
A single inequality constraint has the advantage that it translates into a simpler SOS constraint.
The resulting empty set condition equivalent to ~\bref{eq:clf_emptyset} is then given by:
\begin{align} \label{eq:clf_emptyset_l1}
   \begin{split}
      \Bigl\{ x \in \R^n \mid & \nabla V(x)^T f(x) \geq 0, \nabla V(x)^T G(x) = 0, \\
      &l(x) \neq 0 \Bigr\} = \emptyset.
   \end{split}
\end{align}



According to \tref{thm:positivstellensatz}, by choosing $k=1$, the empty set condition \bref{eq:clf_emptyset_l1} becomes
\begin{align} \label{eq:clf_sos}
   \begin{split}
      - s_1(x) \nabla V(x)^T f(x) - \nabla V(x)^T G(x) p(x) - l(x)^2 \in \Sigma[x],
   \end{split}
\end{align}
where $s_1(x) \in \Sigma[x]$ is an SOS polynomial, and $p(x) = \begin{bmatrix}
   p_1(x) & ... & p_{m}(x)
\end{bmatrix}^T \in (R[x])^m$ is a vector of polynomials.

If $s_1(x)$ is strictly positive w.r.t. $x^*$, the following controller --- in the form of a rational function --- naturally results from the SOS constraint~\bref{eq:clf_sos}:
\begin{align} \label{eq:clf_controller}
      u_{CLF}(x) := p(x)/s_1(x).
\end{align}

\begin{remark} \label{rem:denote_clf}
   Fixing the controller in~\bref{eq:poly_system} to $u_c = u_{CLF}(x)$ results in the closed-loop system~\bref{eq:autonomous_poly_system}.
   Hence, stability is asserted by the existence of a Lyapunov function $V(x)$.
   For the sake of readability, however, we keep denoting $V(x)$ a CLF.
\end{remark}

\begin{lem} \label{lem:clf_stability}
   Given a CLF $V(x)$ for the control system~\bref{eq:poly_system} with $s_1(x)$ in~\bref{eq:clf_sos} strictly positive, the closed-loop system~\bref{eq:autonomous_poly_system} using controller $u(x) = u_{CLF}(x)$ defined in~\bref{eq:clf_controller} is globally asymptotically stable (GAS).
\end{lem}
\begin{pf}
   The polynomial CLF $V(x)$ is strictly positive definite w.r.t. $x^*$.
   Hence, $V(x)$ is radially unbounded. 
   From \bref{eq:clf_sos}, we derive $s_1(x) \nabla V(x)^T f(x) + \nabla V(x)^T G(x) p(x) < 0$ for all $x \neq x^*$. 
   Since $s_1(x)$ is strictly positive w.r.t. $x^*$, the division by $s_1(x)$ results in $\nabla V(x)^T \left ( f(x) + G(x) u_{CLF}(x) \right ) \allowbreak < 0$ for all $x \neq x^*$.
   GAS follows from Theorem 4.2 in \cite{Khalil:02}.
   \hfill$\blacksquare$
\end{pf}

The strict positivity condition on $s_1(x)$ can be enforced by the SOS constraint
\begin{align} \label{eq:sos_s1}
   s_1(x) - \epsilon_{s1} \in \Sigma[x],
\end{align}
for some $\epsilon_{s1} > 0$.

\subsection{Control Barrier Function}

Similar to the stability argument, safety can be asserted with the existence of a scalar function $B_1: \R^n \to \R$.
Specifically, the control system~\bref{eq:poly_system} is safe w.r.t. the set of safe states 
\begin{align} \label{eq:initial_set_1}
   \mathcal X_{s,1} = \Bigl\{ x \in \R^n \mid B_1(x) \leq 0 \Bigr\},
\end{align}
if there exists a control input $u_c \in \R^m$ for every state $x \in \R^n$ such that $\mathcal X_{s,1}$ is forward invariant.
This motivates the following definition of a CBF (cf. \cite{WangMargellos:22}):

\begin{definition} \label{def:cbf}
   Consider a differentiable function $B_1: \R^n \to \R$ such that $\mathcal X_{s,1}$ is non-empty.
   Such a scalar function is called a \underline{Control Barrier Function} (CBF) for the control system~\bref{eq:poly_system} if
   \begin{align} \label{eq:cbf_condition}
      \nabla B_1(x)^T f(x) < 0
   \end{align}
   for all $x \in \Bigl\{ x \in \R^n \mid \nabla B_1(x)^T G(x) = 0, B_1(x) = 0 \Bigr\}$.
\end{definition}
An alternative definition of a CBF (cf. \cite{AmesCoogan:19}) involves a supremum and a class K function.
This definition, however, is not well suited for polynomial optimization since the resulting functions cannot be directly translated to polynomial inequalities.
Similarly to CLF, we restrict our focus on \emph{polynomial} CBFs $B_1(x) \in R[x]$ for the rest of the paper.

Inequality~\bref{eq:cbf_condition} can also be formulated as the empty set condition
\begin{align} \label{eq:barrier_emptyset}
   \begin{split}
       \Bigl\{ x \in \mathbb R^n \mid & \nabla B_1(x)^T f(x) \geq 0, \nabla B_1(x)^T G(x) = 0, \\
       &B_1(x) = 0\Bigr\} = \emptyset.
   \end{split}
\end{align}

\begin{assumption} \label{ass:compact_zero_level_set}
   The set $\bigl\{ x \in \R \mid B_1(x) = 0 \bigr\}$ is compact.
\end{assumption}

According to \tref{thm:putinar_positivstellensatz} and \aref{ass:compact_zero_level_set}, the empty set condition~\bref{eq:barrier_emptyset} is equivalent to the SOS constraint
\begin{align} \label{eq:cbf_sos}
   \begin{split}
       & - s_1(x) \nabla B_1(x)^T f(x) - \nabla B_1(x)^T G(x) p(x) \\
       &\qquad - p_{m+1}(x) B_1(x) \in \Sigma[x],
   \end{split}
\end{align}
where $s_1(x) \in \Sigma[x]$ is an SOS polynomial, $p(x) = \begin{bmatrix}
   p_1(x) & ... & p_{m}(x)
\end{bmatrix}^T \in (R[x])^m$ is a vector of polynomials, and $p_{m+1}(x) \in R[x]$ is a scalar polynomial.

If $s_1(x)$ in~\bref{eq:cbf_sos} is strictly positive w.r.t. $x^*$, there exists a controller defined by the rational function
\begin{align} \label{eq:cbf_controller}
      u_{\text{CBF},1}(x) := p(x)/s_1(x).
\end{align}

\begin{remark}
   As in \rref{rem:denote_clf}, we denote $B(x)$ as CBF for both the control and closed-loop system.
\end{remark}

\begin{lem} \label{lem:forward_invariant}
   Given a CBF $B_1(x)$ for the control system~\bref{eq:poly_system} with $s_1(x)$ in~\bref{eq:cbf_sos} strictly positive, the set $\mathcal X_{s,1}$ in~\bref{eq:initial_set_1} --- if compact (cf.~\aref{ass:compact_zero_level_set}) --- is forward invariant w.r.t. system \bref{eq:autonomous_poly_system} using controller $u(x) = u_{\text{CBF},1}(x)$ defined in~\bref{eq:cbf_controller}.
\end{lem}
\begin{pf}
   A compact zero sublevel set of a scalar polynomial $B_1(x)$ is forward invariant w.r.t. closed-loop system~\bref{eq:autonomous_poly_system} with control $u(x)$, if $\nabla B_1(x)^T \left ( f(x) + G(x) u(x) \right ) < 0$ for all $x \in \lbrace x \in \R^n \mid B_1(x) = 0 \rbrace$ \cite[Nagumo's Theorem 3.1]{Blanchini:99}.
   From~\bref{eq:cbf_sos}, we derive $s_1(x) \nabla B_1(x)^T f(x) + \nabla B_1(x)^T G(x) p(x) < 0$ for all $x \in \lbrace x \in \R^n \mid B_1(x) = 0 \rbrace$.
   Since $s_1(x)$ is strictly positive w.r.t. $x^*$, the division by $s_1(x)$ results in $s_1(x) \nabla B_1(x)^T \left ( f(x) + G(x) u_{\text{CBF},1}(x) \right ) < 0$ for all $x \in \lbrace x \in \R^n \mid B_1(x) = 0 \rbrace$.

   \hfill$\blacksquare$
\end{pf}

\section{SOS Construction of compatible CLF and multiple CBFs} \label{sec:unification}

In this section, we are interested in finding a CLF and multiple CBFs that are compatible with each other.
This is achieved by deriving a new set of SOS constraints that result from restricting the rational controllers \bref{eq:clf_controller} and \bref{eq:cbf_controller} to be identical.

Multiple CBFs provide an additional flexibility to increase the volume of the safe set $\mathcal X_s$, which is defined by the intersection of the zero sublevel sets of all the CBFs.
Intuitively, this makes sense since the solution of a single CBF can be recovered by equating all CBFs.

\subsection{Multiple Control Barrier Functions}

Consider a set of CBFs $\lbrace B_i(x) \rbrace_{i \in [t]}$ for some $t \in \mathbb N_+$, each defining a forward invariant set
\begin{align} \label{eq:initial_set_i}
   \mathcal X_{s,i} = \Bigl\{ x \in \R^n \mid B_i(x) \leq 0 \Bigr\}.
\end{align}
Multiple CBFs are compatible with each other if they have the control-sharing property (cf. \cite{GramBlanch:13}).
\begin{definition} 
   Consider the closed-loop system~\bref{eq:autonomous_poly_system} with feedback controller $u(x)$. A set of CBFs $\lbrace B_i(x) \rbrace_{i \in [t]}$ with associated controllers $u_{\text{CBF},i}(x)$ rendering~\bref{eq:initial_set_i} invariant has the control-sharing property, if $u(x) = u_{\text{CBF},i}(x)$ for all $i \in [t]$.
\end{definition}

The sets of compatible CBFs $\lbrace B_i(x) \rbrace_{i \in [t]}$ form a new forward invariant set $\mathcal X_s$, called safe set, defined as the intersection of the invariant sets from each $B_i(x)$:
\begin{align} \label{eq:initial_set}
   \mathcal X_s := \bigcap_{i \in [t]} \mathcal X_{s,i} = \Bigl\{ x \in \mathbb R^n \mid B_i(x) \leq 0 \text{ } \forall i \in [t] \Bigr\}.
\end{align}

Consider a set of allowable states of the form
\begin{align} \label{eq:allowable_set}
   \mathcal X_a = \Bigl\{ x \in \mathbb R^n \mid w_i(x) \leq 0 \text{ } \forall i \in [t] \Bigr\},
\end{align}
where $w_1(x), ..., w_t \in R[x]$ are polynomials.
Note that for each polynomial $w_i(x)$, that defines the allowable set $\mathcal X_a$, a corresponding CBF can be assigned.
\begin{assumption} \label{ass:compact_allowable_set}
   The sets $\bigl\{ x \in \R \mid w_i(x) \leq 0 \bigr\}_{i \in [t]}$ are compact.
\end{assumption}
To conclude safety, the resulting invariant set $\mathcal X_s$ needs to be contained in the allowable set, i.e. $\mathcal X_s \subseteq \mathcal X_a$.
This can be achieved by restricting the invariant set of each CBF to the zero sublevel set of the corresponding $w_i(x)$ as stated by the following lemma.
\begin{lem} \label{lem:multiple_cbf_safety}
   Consider a set of CBFs $\lbrace B_i(x) \rbrace_{i \in [t]}$ that have the control-sharing property with input $u(x)$.
   If for every $i \in [t]$
   \begin{align} \label{eq:initial_set_subset_wi}
      \mathcal X_{s,i} \subseteq \Bigl\{ x \in \R^n \mid w_i(x) \leq 0 \Bigr\},
   \end{align}
   then the closed-loop system~\bref{eq:autonomous_poly_system} using controller $u(x)$ is safe w.r.t. the safe set $\mathcal X_s$~\bref{eq:initial_set} and the allowable set $\mathcal X_a$~\bref{eq:allowable_set}.
\end{lem}
\begin{pf}
   For a closed-loop system to be safe, we need to show that $\mathcal X_s$ is (i) forward invariant and (ii) $\mathcal X_s \subseteq \mathcal X_a$.
   Note that the sets $\mathcal X_{s,i}$, $i\in[t]$ are compact due to~\bref{eq:initial_set_subset_wi}.
   By \lref{lem:forward_invariant}, every $\mathcal X_{s,i}$, $i\in[t]$, is forward invariant w.r.t. the closed-loop system, and the intersection of forward invariant sets is also forward invariant.
   This proves the condition (i).
   The condition (ii) can be seen from
   \begin{align*}
      \mathcal X_s = \bigcap_{i \in [t]} \mathcal X_{s,i} \subseteq \bigcap_{i \in [t]} \Bigl\{ x \in \R^n \mid w_i(x) \leq 0 \Bigr\} = \mathcal X_a.
   \end{align*}
   \hfill$\blacksquare$
\end{pf}
Finally, we formalize \bref{eq:initial_set_subset_wi} as the empty set condition
\begin{align} \label{eq:initial_subset_allowable}
   \Bigl\{ x \in \mathbb R^n \mid w_i(x) \geq 0, -B_i(x) \geq 0\Bigr\} = \emptyset
\end{align}
for all $i \in [t]$.
Under \tref{thm:putinar_positivstellensatz} and \aref{ass:compact_allowable_set}, condition~\bref{eq:initial_subset_allowable} is equivalent to the following SOS condition:
\begin{align} \label{eq:wi_sos}
   B_i(x) - s_{4,i}(x) w_i(x) \in \Sigma[x]
\end{align}
where $s_{4,i} \in \Sigma[x]$ is an SOS polynomial for all $i \in [t]$.

\subsection{Compatible CLF and CBFs}

The notion of stability and safety can be unified by the existence of both a CLF and multiple CBFs that are compatible with each other via the following control-sharing property.
\begin{definition} \label{def:control_sharing_clf_cbf}
   Consider a closed-loop system~\bref{eq:autonomous_poly_system} with feedback controller $u(x)$. A CLF with associated controller $u_{\text{CLF}}(x)$ stabilizing the system to $x^*$ and set of CBFs $\lbrace B_i(x) \rbrace_{i \in [t]}$ with associated controllers $u_{\text{CBF},i}(x)$ rendering~\bref{eq:initial_set_i} invariant have the control-sharing property, if $u(x) = u_{\text{CLF}}(x)$ and $u(x) = u_{\text{CBF},i}(x)$ for all $i \in [t]$.

\end{definition}

With \dref{def:control_sharing_clf_cbf} in place, we can state our main proposition:
\begin{prop} \label{lem:control_sharing_clf_cbf}
   Given a CLF $V(x)$ and a set of CBFs $\lbrace B_i(x) \rbrace_{i \in [t]}$ that have the control-sharing property with input $u(x)$.
   If condition~\bref{eq:initial_set_subset_wi} holds for every $i \in [t]$, then the closed-loop system~\bref{eq:autonomous_poly_system} with input $u(x)$ is asymptotically stable and safe w.r.t. the allowable set $\mathcal X_a$~\bref{eq:allowable_set}.
\end{prop}
\begin{pf}
   \lref{lem:clf_stability} and \lref{lem:multiple_cbf_safety} prove stability and safety for the closed-loop system using input $u(x)$.

   \hfill$\blacksquare$
\end{pf}

Based on conditions~\bref{eq:clf_sos},~\bref{eq:sos_s1}, ~\bref{eq:cbf_sos} and~\bref{eq:wi_sos}, the stability and safety condition in \lref{lem:control_sharing_clf_cbf} can be restated as the SOS constraints:
\begin{align} \label{eq:sos_condition}
   \begin{array}{lll}
      -\nabla V(x)^T \left( s_1(x) f(x) + G(x) p(x) \right)  - l(x)^2 \\
      \qquad \in \Sigma[x], \\
      -\nabla B_i(x)^T \left( s_1(x) f(x) + G(x) p(x) \right)  \\
      \qquad - p_{m+1,i}(x) B_i(x) \in \Sigma[x], && \forall i \in [t] \\
      B_i(x) - s_{4,i}(x) w_i(x) \in \Sigma[x], && \forall i \in [t] \\
      s_1(x) - \epsilon_{s1} \in \Sigma[x],\\
   \end{array}
\end{align}

where $s_1(x), s_{4,i}(x) \in \Sigma[x]$ are SOS polynomials, $p(x) = \begin{bmatrix}
   p_1(x) & ... & p_{m}(x)
\end{bmatrix}^T \in (R[x])^m$ is a vector of polynomials, and $p_{m+1,i} \in R[x]$ is a scalar polynomial.

The resulting control input $u(x) = p(x) / s_1(x)$ renders the system asymptotically stable on $\mathcal X_a$ and safe w.r.t. the allowable set $\mathcal X_a$ and safe set $\mathcal X_s$.

\section{SOS Program} \label{sec:sos_program}

In this section, we state the SOS program that finds a CLF and multiple CBFs for the control system~\bref{eq:poly_system}.
The SOS program maximizes a surrogate of the volume of $\mathcal X_s$ subject to the SOS condition~\bref{eq:sos_condition}.



\subsection{Cost function} \label{sec:sos_cost}

Ideally, the SOS program would optimize the CLF and CBFs over the volume of the safe set $\mathcal X_s$:
\newline
\begin{align} \label{prg:max_volume}
   \begin{array}{ll}
      \underset{v}{\mbox{minimize}} \quad & -vol \left( \Bigl\{ x \in \mathbb R^n \mid B_i(x, v) \leq 0 \text{ for } i \in [t] \Bigr\} \right) \\
      \mbox{subject to} & \text{SOS constraints~\bref{eq:sos_condition}}
  \end{array}
\end{align}
where $B_i(x, v) = Z_{B,i}(x)^T Q_{B,i}(v) Z_{B,i}(x)$ (cf. \bref{eq:param_poly}).
To make~\bref{prg:max_volume} computationally tractable, the volume of ${\mathcal X}_s$ is surrogated by the traces of the quadratic matrices encoding the $B_i(x)$. The cost then becomes:
\begin{align}
   c_{1}(v) = \sum_{i=1}^t tr(Q_{B,i}(v)).
\end{align}
Furthermore, we have empirically observed that the outcome is improved by restricting the CBF to a predefined center point.
The justification for this lies in the fact that the CBF condition is only valid when $B(x) = 0$ and, therefore, it does not restrict the CBF from taking arbitrary large values when $B(x) \neq 0$.

Let us define a center point $(x_{c,i}, B_{c,i}) \in (\R^n, \R)$ for each CBF $B_i(x),\,i \in [t]$.
Then, we define a cost function by the deviation from that point:
\begin{align}
   c_{2}(v) = \sum_{i=1}^t (B_i(x_{c,i}) - B_{c,i})^2
\end{align}
The cost function of the SOS program is set to
\begin{align} \label{eq:sos_cost}
   c(v) = c_1(v) + c_2(v).
\end{align}

\subsection{SOS Program} \label{sec:final_sos_program}

Given the cost~\bref{eq:sos_cost} and the SOS constraint~\bref{eq:sos_condition}, we summarize the SOS program by:
\begin{align*}
   \hspace{-\leftmargin}
   \begin{array}{lll}
      \mbox{find} & V, B_i, s_1, s_{4,i}, p_{m+1,i} \in R[x], \\
      &p \in (\R[x])^m \\
      \mbox{minimize} & \sum_{i=1}^t tr(Q_{B, i}) + (B_i(x_{c,i}) - B_{c,i})^2 \\
      \mbox{subject to} & -\nabla V^T \left( s_1 f + G  p \right) - l^2 \in \Sigma[x] \\
      &-\nabla B^T_i \left( s_1 f + G  p \right) - p_{m+1,i} B_i \\
      & \qquad \in \Sigma[x] & \forall i \in [t] \\
      & B_i - s_{4,i} w_i \in \Sigma[x] & \forall i \in [t] \\
      & s_1 - \epsilon_{s1} \in \Sigma[x] \\
      & s_{4,i} \in \Sigma[x] & \forall i \in [t] \\
      & V \in \Sigma[x], V(x^*) = 0
   \end{array}
\end{align*}
where $(x_{c,i}, B_{c,i}) \in (\R^n, \R)$ are the predefined center points for each CBF.

\section{Algorithm} \label{sec:algorithm}

The SOS problem in \sref{sec:final_sos_program} cannot directly be converted to an SDP problem due to bilinear terms in its decision variables.
However, an alternating algorithm can find adequate solutions to this non-convex problem.

\subsection{Alternating Algorithm}


Consider an abstract optimization problem with $n_c$ SOS inequality constraints defined by $f_1(v_1, v_2) \in (R[x])^{n_c}$ that are bilinear in the decision variables $v_1 \in \mathbb R^{n_1}$ and $v_2 \in \mathbb R^{n_2}$, and a cost function $f_0(v_1) \in \R$ that only depends on $v_1$:
\begin{align} \label{eq:bilinear_optim}
    \begin{array}{ll}
        \underset{v_1, v_2}{\mbox{minimize}} & f_0(v_1) \\
        \mbox{subject to} & f_1(v_1, v_2) \in (\Sigma[x])^{n_c}
    \end{array}
\end{align}

The optimization problem~\bref{eq:bilinear_optim} cannot be solved directly with an SOS program.
Instead, we propose a method that alternates between searching over one variable while holding fixed the other.

Starting from an iteration $v_1^{(0)}$ and $v_2^{(0)}$, the algorithm at iteration $k$ is defined by:
\begin{itemize}
    \item Step 1: Substitute $v_2=v_2^{(k-1)}$ and solve for $v_1$
    \begin{align} \label{eq:bilinear_largest_margin_step1}
        \begin{array}{ll}
            \underset{v_1}{\mbox{minimize}} & f_0(v_1) \\
            \mbox{subject to} & f_1(v_1, v_2^{(k-1)}) \in (\Sigma[x])^{n_c}
        \end{array}
    \end{align}
    \item Step 2: Substitute $v_1=v_1^{(k)}$ and solve for $\epsilon$ and $v_2$
    \begin{align} \label{eq:bilinear_largest_margin_step2}
        \begin{array}{ll}
            \underset{\epsilon, v_2}{\mbox{minimize}} & \epsilon \\
            \mbox{subject to} & f_1(v_1^{(k)}, v_2) + s(v_1^{(k)}) \epsilon \in (\Sigma[x])^{n_c}
        \end{array}
    \end{align}
\end{itemize}
The algorithm terminates when $|f_0(v_1^{(k)}) - f_0(v_1^{(k-1)})|$ is lower than a given threshold.

\begin{lem}
   If the initial pair $(v_1^{(0)}, v_2^{(0)})$ is feasible for \bref{eq:bilinear_optim}, then feasibility is maintained along the alternating algorithms and the cost $f_0(v_1^{(k)})$ in \bref{eq:bilinear_largest_margin_step1} is non-increasing with iteration $k$.
\end{lem}

\begin{pf}
   First, we show that if Step 1 in iteration $k$ is feasible for \bref{eq:bilinear_largest_margin_step1} then also Step 2 in iteration $k$ must be feasible \bref{eq:bilinear_largest_margin_step2}.
   For feasibility, we only need to find a point $(v_2, \epsilon)$, for which the constraint in Step 2 holds.
   From Step 1, we know that $f_1(v_1^{(k)}, v_2^{(k-1)})$ is an SOS polynomial.
   But then, $(v_2, \epsilon) = (v_2^{(k-1)}, 0)$ is a feasible point of (\ref{eq:bilinear_largest_margin_step2}).
   As a consequence, $\epsilon$ must be equal or smaller than zero, i.e. $\epsilon \leq 0$.

   Next, we show that the solution of iteration $k$ is a feasible point of Step 1 in iteration $k+1$.
   First, note that $-s(v_1^{(k)}) \epsilon$ is an SOS polynomial vector since $\epsilon \leq 0$.
   But then also $f_1 (v_1^{(k)}, v_2^{(k)} )$ is an SOS polynomial.
   Therefore, $(v_1^{(k)}, v_2^{(k)} )$ is a feasible point of Step 1 in iteration $k$.
   Taken together, we showed that $(v_1^{(k)}, v_2^{(k)})$ remain feasible for \bref{eq:bilinear_optim}.
   Therefore , $f_0(v_1^{(k+1)})$ is either equal or smaller than $f_0(v_1^{(k)})$.

   \hfill$\blacksquare$
\end{pf}

\subsection{Main algorithm}  \label{sec:main_algorithm}

The main algorithm improves feasibility by introducing an operating region $\mathcal X_{op}$.
Under \aref{ass:compact_allowable_set}, there exists a scalar polynomial $r(x)$ whose zero sublevel set is compact and strictly contains the allowable set $\mathcal X_a$, i.e. 
\begin{align} \label{eq:operating_region}
   \mathcal X_s \subseteq \mathcal X_a \subseteq{\mathcal X}_{op} := \Bigl\{ x \in \mathbb R^n \mid r(x) \leq 0 \Bigr\}.
\end{align}
Note that all state trajectories of interest are contained in the operating region $\mathcal X_{op}$.
We will use this fact to slightly alter the empty set conditions and the corresponding SOS constraints developed in the previous chapters.

Given an allowable set $\mathcal X_a$ defined by $w_i(x)$ for $i \in [t]$~\bref{eq:allowable_set}, positive definite SOS polynomials $l(x) \in \Sigma[x]$, a positive scalar $\epsilon_{s1} > 0$, center points $(x_{c,i}, B_{c,i}) \in (\R^n, \R)$ for each $i \in [t]$, a region of operation defined by $r(x)$ \bref{eq:operating_region} and initial values $s_1(x)$, $p(x)$, and $p_{m+1,i}(x)$ for $i \in [t]$, the main SOS algorithm is summarized below.

\begin{itemize}
   \item Step 1: Given a controller $u = p / s_1$ find a CLF and multiple CBFs.
   \begin{align*}
       \hspace{-\leftmargin}
       \begin{array}{lll}
           \mbox{find} & V, B_i, s_2, s_{3,i}, s_{4,i} \in R[x] \\
           \mbox{hold fixed} \quad & s_1, p_{m+1,i} \in R[x], p \in (\R[x])^m \\
           \mbox{minimize} & \sum_{i=1}^t tr(Q_{B, i}) + (B_i(x_{c,i}) - B_{c,i})^2 \\
           \mbox{subject to} & -\nabla V^T \left( s_1 f + G p \right)  - l^2 + s_2 r \in \Sigma[x] \\
           &-\nabla B^T_i \left( s_1 f + G p \right) - p_{m+1,i} B_i \\
           &\qquad+ s_{3,i} r \in \Sigma[x] & \forall i \in [t] \\
           & B_i - s_{4,i} w_i \in \Sigma[x] & \forall i \in [t] \\
           & s_2 \in \Sigma[x] \\
           & s_{3,i} \in \Sigma[x] & \forall i \in [t] \\
           & s_{4,i} \in \Sigma[x] & \forall i \in [t] \\
           & V \in \Sigma[x], V(x^*) = 0 \\
           \\
       \end{array}
   \end{align*}

   \item Step 2: Given a CLF and multiple CBFs find a controller $u = p / s_1$.
   \begin{align*}
       \hspace{-\leftmargin}
       \begin{array}{lll}
           \mbox{find} & \epsilon_i \in \R, s_1, s_2, s_{3,i}, p_{m+1,i} \in R[x], \\
           & p \in (\R[x])^m \\
           \mbox{hold fixed} \quad & V, B_i \in R[x] \\
           \mbox{minimize} & \sum_{i=1}^t \epsilon_i \\
           \mbox{subject to} & -\nabla V^T \left( s_1 f + G p \right) - l^2 + s_2 r \in \Sigma[x] \\
           & -\nabla B^T_i \left( s_1 f + G p \right) - p_{m+1,i} B_i \\ 
           &\qquad + s_{3,i} r + \epsilon_i \in \Sigma[x] & \forall i \in [t] \\
           & s_1 - \epsilon_{s1} \in \Sigma[x] \\
           & s_2 \in \Sigma[x] \\
           & s_{3,i} \in \Sigma[x] & \forall i \in [t] \\
       \end{array}
   \end{align*}
\end{itemize}

\section{Simulation} \label{sec:simulation}

The alternating algorithm from~\sref{sec:main_algorithm} is implemented for a three-dimensional vector field realization of a dc/ac power converter model.
A feedback controller that avoids unsafe states is of paramount importance for this application since all electrical variables (e.g. voltages and currents) need to be constrained at all times.
With the dc voltage and ac currents as states, the control system~\bref{eq:poly_system} is defined by 
\begin{align}
   f(x) = \begin{bmatrix}
      -0.05 x_1 - 57.9 x_2 + 0.00919 x_3 \\
      1710 x_1 + 314 x_3 \\
      -0.271 x_1 - 314 x_2
   \end{bmatrix}
\end{align}
and
\begin{align}
   G(x) = \begin{bmatrix}
      0.05 - 57.9 x_2 & -57.9 x_3 \\
      1710 + 1710 x_1 & 0 \\
      0 & 1710 + 1710 x_1
   \end{bmatrix}.
\end{align}
The allowable set $\mathcal X_a$ is defined by two polynomials $w_1(x)$ and $w_2(x)$ encoding the state constraints for current and voltage respectively:
\begin{align}
   \begin{split}
      w_1(x) &= (x_1 + 0.3)^2 + (x_2 / 20)^2 + (x_3/20)^2 - 0.5^2 \\
      w_2(x) &= (x_1/20)^2 + x_2^2 + x_3^2 - 1.2^2.
   \end{split}
\end{align}

An operating region $\mathcal X_{op}$ in~\bref{eq:operating_region} can be parameterized by:
\begin{align}
   r(x) = (x_1/0.8)^2 + (x_2/1.2)^2 + (x_3/1.2)^2 - 1.8.
\end{align}

A $r(x)$ that tighter fits the allowable set results in a faster convergence of the alternating algorithm.
For illustration purposes, however, we selected a suboptimal operation region.
The center points of the CBFs are given by:
$x_{c,i} = \begin{bmatrix}
   -0.3 & 0 & 0
\end{bmatrix}^T$
 and $B_{c,i} = -10$ for $i \in [2]$.
The degrees of the polynomials involved in the algorithm are summarized as follows:
\begin{center}
   \begin{tabular}{||c@{\hskip 0.04in} c@{\hskip 0.04in} c@{\hskip 0.04in} c@{\hskip 0.04in} c@{\hskip 0.04in} c@{\hskip 0.04in} c@{\hskip 0.04in} c||} \hline
    $V(x)$ & $B_i(x)$ & $s_1(x)$ & $s_2(x)$ & $s_{3,i}(x)$ & $s_{4,i}(x)$ & $p(x)$ & $p_{m+1,i}(x)$ \\ \hline\hline
    4 & 4 & 2 & 6 & 6 & 3 & 4 & 2 \\ \hline
   \end{tabular}
\end{center}
For $V(x)$ and $B_i(x)$, the odd degree coefficients are set to zero:
\begin{itemize}
   \item $deg \left( V(x) \right) = (2, 4)$
   \item $deg \left( B_1(x - x_{c,1}) \right) = (0, 2, 4)$
   \item $deg \left( B_2(x - x_{c,2}) \right) = (0, 2, 4)$.
\end{itemize}

\begin{figure}[htbp] \label{fig:evolution_safe_set}
   \centerline{\includegraphics[width=0.5\textwidth]{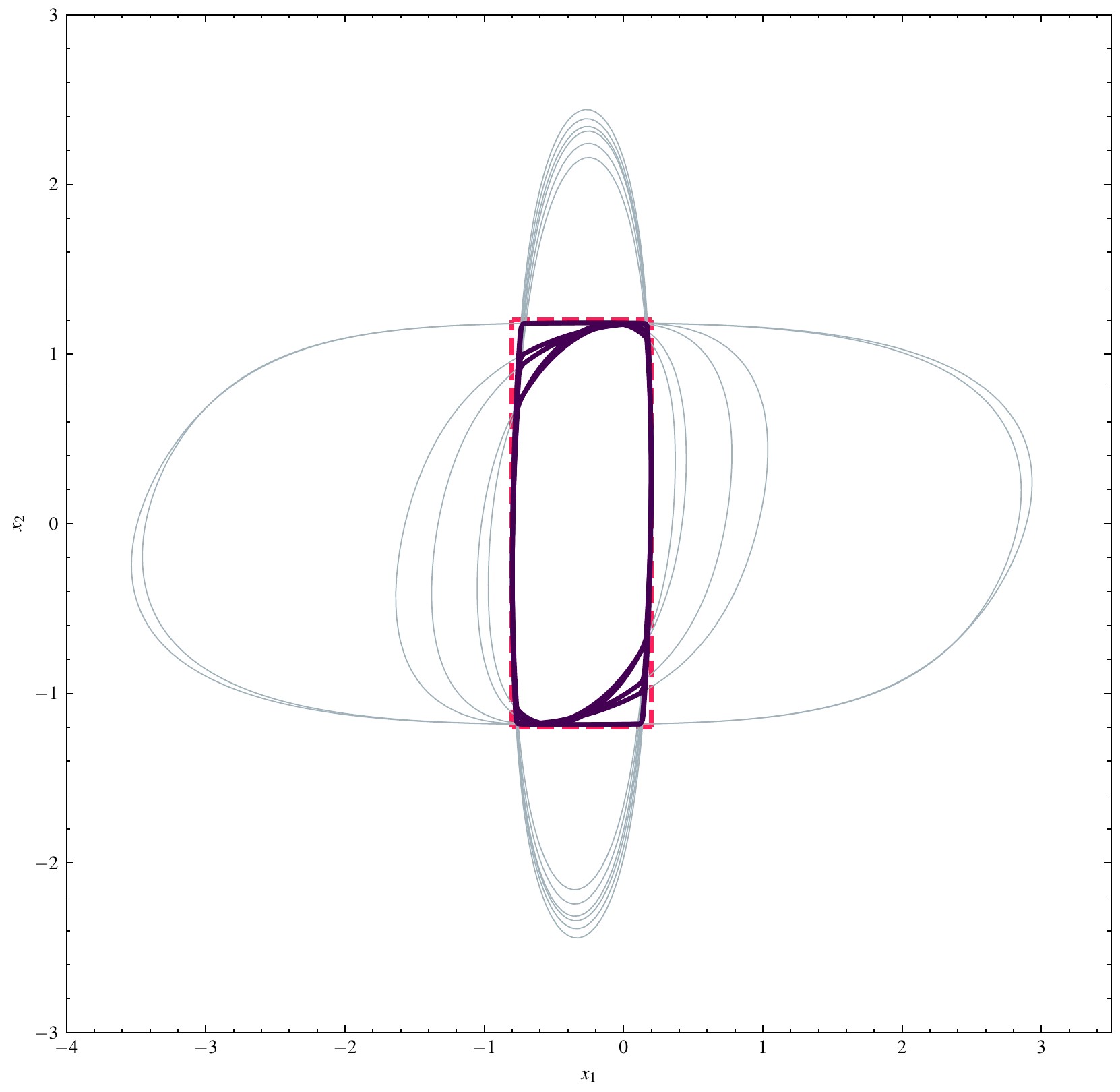}}
   \caption{shows the evolution of the sequence of safe sets $\mathcal X_s$ projected to the $x_1$ and $x_2$ coordinates over 6 iterations. The boundary of the allowable set $\mathcal X_a$ is shown in red. The sequence of zero-level sets of the CBF $B_1(x)$ and $B_2(x)$ are indicated in gray.}
\end{figure}

\fref{fig:evolution_safe_set} illustrates the sequence of safe sets calculated at each iteration of the algorithm.
The computed CLF and CBFs not only prove that there exists a stabilizing feedback controller but also guarantee that such controllers renders the system safe w.r.t. a safe set $\mathcal X_s$.
By using two CBFs, the algorithm terminates with a safe set that tightly fits the allowable set $\mathcal X_a$.

The number of variables and constraints involved in the resulting SDP are summarized as follows:
\begin{center}
   \begin{tabular}{||ccc||} 
      \hline
      & Step 1 & Step 2 \\ \hline\hline
      Variables & 337 & 372 \\ \hline
      Constraints & 5188 & 4892 \\
      \hline
   \end{tabular}
\end{center}
The algorithm was completed in 90 seconds.
For each iteration, the SDP algorithm took on average 33 inner iterations to solve step 1 and on average 21 inner iterations to solve step 2.

\section{Conclusions and future work} \label{sec:conclusions}

This paper presented a framework to combine stability and safety conditions using CLF and multiple CBFs.
By employing two versions of the Positivstellensatz, we synthesized a controller used to prove compatibility between CLF and CBFs.
We then formalized SOS constraints that encode compatible CLF and CBF conditions.
Finally, we proposed an algorithm that solves the resulting SOS program by iteratively solving two SDPs.

For future work, we plan to study the computational complexity of our algorithm and also to explore a unified framework, that proves stability and safety with weaker CLF and CBF conditions. 
Moreover, we intend to incorporate noise-robustness into our SOS formulation, as demonstrated by \cite{KangChen:23}.

\bibliography{ifacconf}

\end{document}